\documentclass[11pt]{amsart}
\usepackage{amsfonts,amssymb,eucal}
\usepackage[all]{xy}
\begin{document}

\newcommand{\C}{{\mathbb C}}
\newcommand{\F}{{\mathbb F}}
\newcommand{\R}{{\mathbb R}}
\newcommand{\Q}{{\mathbb Q}}
\newcommand{\Z}{{\mathbb Z}}
\newcommand{\hh}{{\mathbb H}}
\newcommand{\Oh}{{\rm O}}
\newcommand{\Spin}{{\rm Spin}}
\newcommand{\SO}{{\rm SO}}
\newcommand{\SU}{{\rm SU}}
\newcommand{\U}{{\rm U}}
\newcommand{\T}{{\mathbb T}}
\newcommand{\ad}{{\rm ad}}
\newcommand{\op}{{\rm op}}

\newcommand{\bC}{{\bf C}}
\newcommand{\bD}{{\bf D}}
\newcommand{\bH}{{\bf H}}
\newcommand{\cC}{{\mathcal C}}
\newcommand{\rL}{{\rm L}}

\newcommand{\cobord}{{\sf Cobord}}
\newcommand{\diff}{{\rm Diff}}
\newcommand{\Metrics}{{\rm Metrics}}
\newcommand{\pt}{{\rm pt}}
\newcommand{\Spec}{{\rm Spec}}
\newcommand{\vect}{{\rm Vect}}
\newcommand{\MT}{{\rm MT}}
\newcommand{\Pic}{{\rm Pic}}
\newcommand{\Aut}{{\rm Aut}}
\newcommand{\vv}{{\rm v}}

\title{Spin Cobordism Categories in Low Dimensions}
\author{Nitu Kitchloo and Jack Morava}

\thanks{Both authors were supported by the NSF}
\subjclass{}
\date {9 November 2006}

\begin{abstract}
{The Madsen-Tillmann spectra defined by categories of
three- and four-dimensional Spin manifolds have a very rich
algebraic structure, whose surface is scratched here.}
\end{abstract}

\maketitle

\begin{center}
\it For Michael Atiyah, in deep gratitude.
\end{center}

\section{Cobordism categories}

\noindent
{\bf 1.1} Many variations and generalizations are possible, but to begin, consider the
topological two-category $D\cobord$ whose objects are oriented smooth closed $d$-manifolds
($D = d+1$), with the topological category $D\cobord(V,V')$ of morphisms having as objects,
$D$-dimensional cobordisms
\[
W : V \to V'
\]
from $V$ to $V'$; for our purposes this will mean an identification $\partial W \cong V_\op
\coprod V'$, extended to a small neighborhood of the boundary. The two-morphisms will be
orientation-preserving diffeomorphisms of such cobordisms, which equal the identity
near the boundary. The composition functor
\[
D\cobord(V,V') \times D\cobord(V',V'') \to D\cobord(V,V'')
\]
is defined by glueing outgoing to incoming boundaries. \bigskip

\noindent
A topological category $\bC$ is a kind of simplicial space, and so has a geometric realization
$|\bC|$; for example,
\[
|D\cobord(V,V')| = \coprod_{[W: V \to V']} B\diff_0^+(W)
\]
is the union, indexed by diffeomorphism classes of cobordisms $W$ from $V$ to $V'$, of the
classifying spaces of the groups of orientation-preserving diffeomorphisms of $W$ which equal
the identity near the boundary. We'll write $D|\cobord|$ for the topological category with closed
$d$-manifolds as objects, and the classifying spaces above as morphism objects; it is symmetric
monoidal (under disjoint union). \bigskip

\noindent
Such categories have an impressive history [3, 16, 17, 21] in topology and physics. This paper applies
the recent breakthroughs of [9] which (in great generality) identify the classifying
spectra of such categories. The formalism of Galatius, Madsen, Tillmann, and Weiss frames these
categories somewhat differently: they work with a category $\cC_D$ of manifolds embedded in a
high-dimensional Euclidean background, but the description used above is equivalent, and is
convenient in physics. \bigskip

\noindent
{\bf 1.2} A topological transformation group $G \times X \to X$ has an associated
homotopy-to-geometric quotient map
\[
X//G := EG \times_G X \to \pt \times_G X = X/G
\]
which defines a kind of resolution
\[
B\diff \sim E\diff \times_\diff \Metrics \to \pt \times_\diff  \Metrics
\]
of the moduli space of Riemannian metrics [8] on a manifold. For a closed manifold the action
of the diffeomorphism group on the space of metrics is proper; for surfaces of genus $> 1$, for
example, its isotropy groups are not just compact but finite, making the map a rational homology
equivalence. \bigskip

\noindent
This resolution defines a monoidal functor
\[
D|\cobord| \to {\rm Gravity}_D
\]
to a topological category with moduli spaces of metrics as its morphism objects. The
Einstein-Hilbert functional
\[
g \mapsto \int_W R(g) \; d {\rm vol}_g : \Metrics/\diff \to \R
\]
is a natural candidate for a Morse function on these objects, so this category models interesting
aspects of (Euclidean) general relativity. Witten has suggested that backgrounding the choice
of Morse function leads to more general models in which topology change can be treated quite
naturally. \bigskip

\noindent
{\bf 1.3.1} This paper is concerned with the cobordism category defined by four-dimensional
Spin manifolds. The classifying space $|\bC|$ of a {\bf symmetric monoidal} topological
category $\bC$ is a kind of abelian monoid, or better: a $\Gamma$-space in the sense of [20].
Its group completion
\[
|\bC|^+ := \Omega B |\bC|
\]
is an infinite loop-space, which is characterized by its associated stable spectrum. \bigskip

\noindent
GMTW identify $|D|\cobord||^+$ as the infinite loopspace associated to a twisted desuspension
\[
MT\SO(D) := B\SO(D)^{-\bD}
\]
of the classifying space for the orthogonal group, where $\bD$ is the vector bundle associated to
the basic representation of $\SO(D)$ on $\R^D$. More generally, a pullback diagram
\[
\xymatrix{
G(d) \ar[d] \ar[r] & \SO(d) \ar[d] \\
G(D) \ar[r] & \SO(D) }
\]
of groups and homomorphisms defines a topological category $G|\cobord|$ of manifolds with
$G(d)$-structures on their tangent bundles, up to cobordism through manifolds with
$G(D)$-structures; and the techniques of [9] identify its associated spectrum as
\[
MTG(D) : = BG(D)^{-\bD}
\]
where $\bD$ is now the $D$-dimensional representation of $G(D)$ pulled back from the basic
representation of the orthogonal group. A further generalization identifies the classifying
spectrum for the category of $G$-manifolds mapped to some parameter space $X$ as
\[
X_+ \wedge BG(D)^{-\bD} \;.
\] \bigskip

\noindent
{\bf 1.3.2} When $d=1$, for example, we get the desuspension (alternately denoted $\C
P^\infty_{-1}$) of  $B\SO(2) = \C P^\infty$ by the tautological line bundle. Its homology is free
of rank one in even dimensions $\geq -2$, and is otherwise zero. Pinching off its
(-2)-dimensional cell defines a cofibration
\[
S^{-2} \to \C P^\infty_{-1} \to \C P^\infty_+
\]
of spectra, with an associated fibration
\[
\Omega^\infty S^{-2} \to |2|\cobord||^+ \to Q(\C P^\infty) \times Q(S^0)
\]
of loopspaces. Since $\Omega^\infty S^{-2}$ has torsion homotopy, the rational homology of
$|2|\cobord||^+$ is a free bicommutative Hopf algebra generated by $\tilde{H}_*\C P_\infty$ (ie, by
the Miller-Morita-Mumford classes $\kappa_i, \; i \geq 1$), extended by degree zero classes
$\kappa_0^{\pm n}$ coming from the rational cohomology of $Q(S^0)$ \bigskip

\noindent
{\bf 1.3.3} The covariant functor [12 Ch III]
\[
X \mapsto \Spec \; H^\pm(\Omega^\infty_0 X,\F) := \tilde{\bH}_\pm(X,\F)
\]
(from (Spectra) to unipotent commutative supergroup-schemes over the field $\F$) is a
homotopy-theoretic analog of the `big' quantum cohomology studied in some contexts in
physics. For example, the infinite loopspace associated to the suspension spectrum
$\Sigma^\infty X$ defined by a connected pointed space splits stably as
\[
\Omega^\infty \Sigma^\infty X \sim \coprod_{n \geq 0} E\Sigma_n \wedge_{\Sigma_n}
X^{\wedge n}
\]
so its rational cohomology is the symmetric algebra on the reduced cohomology of $X$. In this
case $\tilde{\bH}_\pm(X,\Q)$ can be identified with the affine (super)groupscheme
which represents the functor
\[
( \Q - {\rm algebras}) \ni A \mapsto \tilde{H}_\pm(X;A) \in (\Q - \vect) \;.
\]
For a general connected spectrum $X$, $H^\pm(\Omega^\infty X,\Q)$ is the universal enveloping
Hopf algebra associated to the (super-commutative) Lie algebra $\pi_\pm(X) \otimes \Q$ of
primitives. The category of such affine groupschemes is closed and symmetric monoidal, with a
product $\boxtimes$ which is not very familiar [12 Ch II]; over $\Q$, it corresponds to the
graded tensor product of spaces of primitives. \bigskip

\noindent
If the loop-space associated to a spectrum is not connected, let $\bH^0(X,\F)$ be the
groupscheme represented by the group ring $\F[\pi_0(X)]$, and let $\bH_0(X,\F)$ be the
spectrum of the ring of finitely-supported $\F$-valued functions on $\pi_0(X)$; then we can
define
\[
\bH_\pm(X,\F) = \bH_0(X,\F) \times \tilde{\bH}_\pm(X,\F)
\]
(and similarly, for cohomology). For example,
\[
\bH_\pm(S^2 \times MT\SO(2),\Q)
\]
is an analog of the big quantum cohomology related to the Toda lattice [10]. \bigskip

\section{Low-dimensional Spin cobordisms}\bigskip

\noindent
{\bf 2.0} The action
\[
u,q \mapsto uqu^{-1} : \SU(2) \times \hh \to \hh
\]
of the group $\SU(2) = \{ u \in \hh \; | \; |u| = 1 \}$ of unit quaternions leaves the subspace $\R
\subset \hh$ invariant, defining a double cover
\[
\rho : \SU(2) \to \SO(3)
\]
of the rotation group of the subspace orthogonal to it, identifying $\SU(2)$ with $\Spin(3)$.
Similarly, the action
\[
u_L, u_R, q \mapsto u_L q u_R^{-1} : (\SU(2) \times \SU(2)) \times \hh \to \hh
\]
factors through the double cover
\[
\SU(2) \times SU(2) = \Spin(4) \to \SO(4) \;.
\]
It is easy to see that the diagram
\[
\xymatrix{
\SU(2) \ar[r]^\rho \ar[d]^\Delta & \SO(3) \ar[d] \\
\SU(2) \times \SU(2) \ar[r] & \SO(4) }
\]
is a pullback; following \S 1.3.1, this defines the cobordism category of Spin
three-manifolds up to four-dimensional Spin cobordism. Similarly, the $D=3$ Spin cobordism
category is defined by three-dimensional Spin cobordisms between two-dimensional Spin
manifolds: in Riemann surface terms [2], the latter structure amounts to a choice of
square root for the canonical complex line bundle. [Complex Spin structures are very interesting
[22], but they won't be considered here.] \bigskip

\noindent
We'll write $\hh_\ad - 1$ for the three-dimensional representation $\rho$, and $\hh \otimes \hh_\op$
for the four-dimensional Spin representation, as in [13 \S 1.4, 15 \S 1]; then
\[
MT\Spin(3) \sim \Sigma B\SU(2)^{-\hh_\ad}
\]
and
\[
MT\Spin(4) \sim B(\SU(2) \times \SU(2))^{-\hh \otimes \hh_\op} \;.
\] \bigskip

\noindent
These spectra are very nice, with torsion-free integral homology concentrated in degrees
$\equiv - 1$ (resp. $0$) mod four, but they are nontrivial in negative dimensions,
starting in degree $-3$ (resp. $-4$). It will simplify notation below to introduce their
connective suspensions
\[
\MT(3) := \Sigma^3 MT\Spin(3)
\]
and
\[
\MT(4) := \Sigma^4 MT\Spin(4) \;.
\] \bigskip

\noindent
{\bf 2.1} The representation $\hh \otimes \hh_\op$ restricts to $\hh_\ad$ along the diagonal
embedding of $\SU(2)$ in $\SU(2) \times \SU(2)$. Since Thom spaces (and spectra) behave
nicely under pullback, this defines a morphism
\[
\Delta_\natural : \MT(3) \to \MT(4) \;.
\]
The main result of this note asserts that (at least, up to cohomology) this map makes $\MT(3)$
a kind of cocommutative and coassociative coalgebra spectrum. \bigskip

\noindent
{\sc Proposition:} The integral cohomology $H^*\MT(3)$ can be identified with $H^*B\SU(2)$
{\bf as an algebra}, consistent with a splitting
\[
\Psi^* : H^*\MT(3) \otimes H^*\MT(3) \cong H^*\MT(4)
\]
which identifies $\Delta^*_\natural : H^*\MT(4) \to H^*\MT(3)$ with the multiplication map.
\bigskip

\noindent
{\bf Proof:} If $X$ is a compact connected space, then any $[V] \in \tilde{K}\Oh(X)$ is
stably equivalent to a vector bundle $V$ over $X$, of dimension $\vv \gg 0$, and Atiyah's
Thom spectrum
\[
X^{[V]} := \Sigma^{-\vv} X^V
\]
is well-defined up to homotopy. If $[V]$ is orientable (eg if $w_1(V) = 0$, in the case of integral
homology), there is a Thom isomorphism
\[
\Phi_V : H^*X \to H^*X^{[V]} \;.
\]
Taking a limit over finite subcomplexes extends such constructions to nice spaces like
$B\SU(2)$. \bigskip

\noindent
With this notation, we have a commutative diagram

\[
\xymatrix{
H^*MT(3) \otimes H^*\MT(3) \ar[r]^{\Psi^*} & H^*\MT(4) \ar[r]^{\Delta^*_\natural} & H^*\MT(3) \\
H^*B\SU(2)_+ \otimes H^*B\SU(2)_+ \ar[u]^{\Phi_{-\hh_\ad} \otimes \Phi_{-\hh_\ad}} \ar[r]^\cong
& H^*B(\SU(2) \times \SU(2))_+ \ar[u]^{\Phi_{-\hh \otimes \hh_\op}} \ar[r]^{\Delta^*} & H^*B\SU(2)+
\ar[u]^{\Phi_{-\hh_\ad}} }
\] \bigskip

\noindent
with the composition $\Delta^*_\natural \circ \Psi^*$ defining the multiplicative structure.
\bigskip

\noindent
Verification of associativity amounts to unwinding the collection of Thom isomorphisms which
reduce the commutativity (after taking cohomology) of the diagram

\[
\xymatrix{
\MT(3) \ar[r]^{\Delta_\natural} & \MT(4) \ar[r]^<\Psi & \MT(3) \wedge \MT(3) \ar[r]^{1 \wedge
\Delta_\natural} \ar[d]^{\Delta^\natural \wedge 1} & \MT(3) \wedge \MT(4) \ar[d]^{1 \wedge
\Psi}\\
{} & {} & \MT(4) \wedge MT(3) \ar[r]^<{\Psi \wedge 1} & \MT(3) \wedge \MT(3) \wedge
\MT(3) }
\]

\noindent
to a similar diagram expressing the associativity of the usual multiplication on $H^*B\SU(2)_+$.
Commutativity is a consequence of $\Delta_\natural$ being essentially a diagonal, and the
unit is the composition
\[
\xymatrix{
H^*B\SU(2)^{-[\hh_\ad]} \ar[r]^{\Phi^{-1}_{\hh_\ad}} & H^*B\SU(2)_+ \ar[r] & H^*S^0 }
\]
defined by the inclusion of a basepoint into $B\SU(2). \; \; \; \; \; \Box$ \bigskip

\noindent
{\bf 2.2.1} The result above can also be paraphrased in terms of a ring structure on homotopy
quantum cohomology, but because
\[
\Omega^\infty \MT(3) \sim \Z \times \Omega^\infty_0 \MT(3)
\]
is not connected, this requires some discussion. According to \S 1.3.1, $\MT(3)$ is the
cobordism spectrum of Spin three-manifolds mapped to the three-sphere. The extra data defined
by such a map is (at least, after tensoring with $\Q$) very close to a framing (in the sense of
[4, 15 \S 2.1]) of a three-dimensional Spin cobordism. \bigskip

\noindent
Similarly,
\[
\Omega^\infty MT\Spin(4) \sim \Z^2 \times \Omega^\infty MT\Spin(4)
\]
with the classes of a K3 surface and the quaternionic projective plane as natural geometric
generators for $\pi_0$ [11]. The Euler characteristic $\chi$ and the signature $\sigma$ are a
basis for the linear functionals on this group, at least over $\Z[1/2]$, and if $\chi^*,\sigma^*$ denote
the dual basis elements, then

\[
\left[\begin{array}{c}
       \hh P_2 \\ K3
      \end{array}\right] =
\left[\begin{array}{cc}
        2 & 0 \\
        6 & 16
       \end{array}\right]
\left[\begin{array}{c}
        \chi^* \\ \sigma^*
        \end{array}\right] \;.
\]
\bigskip

\noindent
Desuspending the isomorphisms in \S 2.1 yields a splitting
\[
\bH^\pm(MT\Spin(4),\Q) \cong \boxtimes^2 (\bH^{\pm}_{ KM}(\Sigma^{-1} MT\Spin(3),\Q))
\]
of homotopy-theoretic {quantum \bf co}homology: where the subscript on the right indicates an extension
of $\tilde{\bH}^\pm(\Sigma^{-1}MT\Spin(3),\Q)$ by the multiplicative group
(represented by the group ring of Kirby-Melvin framings). \bigskip

\noindent
{\bf 2.2.2} The existence of a multiplication on $H^*\MT(3)$ raises the possibility of the
existence of a so-called Hopf algebroid structure on $(H^*\MT(3),H^*\MT(4))$. In fact, two
three-dimensional cobordisms mapped to the three-sphere define a fiber product
\[
\xymatrix{
W_0 \times_{S^3} W_1 \ar@{.>}[d] \ar@{.>}[r] & W_0 \times W_1 \ar[d] \\
S^3 \ar[r]^\Delta & S^3 \times S^3 }
\]
which is generically another such; but whether this can be used to define a geometric
product on $\MT(3)$ involves subtle questions about framings. \bigskip

\noindent
Note that the symplectic pairing (and the associated duality) on the Tate cohomology $t_\T H\Z$
studied in connection with $MT\SO(2)$ in [18] has a very nice analog on $t_{\SU(2)}H\Z$. \bigskip

\noindent
{\bf 2.2.3} The spectrum
\[
MT\rL(4) := B\SU(2)^{-\hh}
\]
(defined by the obvious action $\sigma$ of $\SU(2)$ by left multiplication on $\hh$) is the
Madsen-Tillmann spectrum of the category of three-dimensional Spin manifolds, up to
cobordism through four-manifolds with an `almost hyperHermitian' structure (in the sense of
[7]). Its cobordisms are essentially four-manifolds with ${\rm Sl}_2(\C)$ (ie, Lorentzian Spin)
structures; its cohomology is concentrated in even dimensions, but it has no very obvious
multiplicative structure. It would be interesting to understand better the relations between
this spectrum and $MT\Spin(3)$, which has cohomology concentrated in odd degrees: it is tempting
to think of $MT\rL(4)$ as some kind of bosonization of $MT\Spin(3)$. \bigskip

\noindent
Behind this lie broader questions about Atiyah-twistings of spectra: the isomorphism
\[
\xymatrix{
H^*B\SU(2)^{-\hh} \ar[r]^{\Phi^{-1}_{-\hh}} & H^*B\SU(2)_+ \ar[r]^{\Phi_{-\hh_\ad}} &
H^*B\SU(2)^{-\hh_\ad} }
\]
does not respect Steenrod operations. In general, a vector bundle $V \to X$
which is oriented with respect to a reasonable multiplicative cohomology theory $E^*$ defines
a rank one projective $E^*(X)$-module $E^*(X^V)$, and thus an element of the Picard group of $E^*(X)$.
These groups tend to be trivial, but their equivariant analogs (with respect to the
cohomology automorphisms of $E$) can be more interesting. \bigskip

\noindent
The spherical fibration associated to $V$ defines a natural invariant
\[
\Pic_{\Aut(E)}(E^*(X)) \to H^1(\Aut(E),(1 + \tilde{E}^*(X))^\times)
\]
which can be pulled back to universal examples involving the $J$-groups
of classifying spaces [6]. Techniques developed for the circle group [15]
seem promising for $\SU(2)$ as well. \bigskip

\noindent
{\bf 2.2.4} Since this paper was submitted, J. Lurie's important work on topological
field theories has become available. We close by drawing attention to some
applications of his ideas to the subject of this paper. \bigskip

\noindent
Lurie's Theorem 2.5.10 [16] identifies the space of infinite-loop maps from $|G|\cobord||^+$
to an infinite loopspace $X = \{X_n\}$ as the homotopy fixed-point spectrum $X^{hG}$ associated
to an action of $G$ on $X$ via the natural action of $\SO(D)$ on suspension coordinates of the stabilization 
of $\Sigma^{D} X_{n-D}$. [This action is closely related to the constructions in the preceding paragraph.] 
The $n$th space of the fixed-point spectrum $X^{hG}$ is equivalent to the space of maps from the Thom 
spectrum $MTG(D)$ (in the notation of \S 1.3.1) to $X_{n-D}$. \bigskip

\noindent
The infinite loopspace $B\otimes$ associated to the monoidal category of real vector spaces under
tensor product is an interesting example. A monoidal functor from $G|\cobord|$ to (Vect,$\otimes$)
is a generalization of a topological quantum field theory in Atiyah's sense, and it defines an
infinite-loop map from $|G|\cobord||^+$ to $B\otimes$, and hence an element of $k_\otimes^{-D}(MTG(D))$.
These groups are accessible via the Atiyah-Segal exponential [5, 19].

\newpage

\bibliographystyle{amsplain}

\bigskip

\noindent 
University of California at San Diego, {\tt nitu@math.ucsd.edu}

\medskip

\noindent Johns Hopkins University, {\tt jack@math.jhu.edu}

\end{document}